\documentclass[preprint,10pt]{elsarticle}
\usepackage{bbding}
\usepackage{indentfirst}
\usepackage{graphicx}
\usepackage{epstopdf}
\usepackage{amsthm,amsmath}
\usepackage{mathrsfs}
\biboptions{numbers,sort&compress}

\usepackage{amssymb}
\UseRawInputEncoding

\begin{document}

\begin{frontmatter}



\title{Asymptotic and Finite-time Cluster Synchronization of Neural Networks via Two Different Controllers}

\author[rvt]{Juan Cao}
\ead{caojuan@nuaa.edu.cn}
\author[rvt]{Fengli Ren\corref{cor}}
\ead{flren@nuaa.edu.cn}
\author[els]{Dacheng Zhou}
\ead{zhoudacheng6@163.com}


\cortext[cor]{Corresponding author}


\address[rvt]{Department of Mathematics, Nanjing University of Aeronautics and Astronautics, Nanjing 210016, P.R. China}
\address[els]{School of Mechatronics Engineering and Automation, Shanghai University, Shanghai 200444, P.R. China}

\begin{abstract}
In this paper, by using pinning impulse controller and hybrid controller respectively, the research difficulties of asymptotic synchronization and finite time cluster synchronization of time-varying delayed neural networks are studied. On the ground of Lyapunov stability theorem and Lyapunov-Razumikhin method, a novel sufficient criterion on asymptotic cluster synchronization of time-varying delayed neural networks is obtained. Utilizing Finite time stability theorem and hybrid control technology, a sufficient criterion on finite-time cluster synchronization is also obtained. In order to deal with time-varying delay and save control cost, pinning pulse control is introduced to promote the realization of asymptotic cluster synchronization. Following the idea of pinning control scheme, we design a progressive hybrid control to promote the realization of finite time cluster synchronization. Finally, an example is given to illustrate the theoretical results.
\end{abstract}


\end{frontmatter}


\section{Introduction}
\par Synchronization in complex networks refers to a evolutionary process where two or more chaotic systems improve their behavior to achieve a common phenomenon under certain coupling conditions or under certain forces. It's an important topic that have been captured by many experts and scholars from various fields of research \cite{s1,s2,s3,s4,s5,s6}, such as physical science, natural science, mathematics, communication, and engineering. Because of that, different types of synchronization protocols have been studied in detail. For example, projective synchronization \cite{s7,s8,s9}, cluster synchronization \cite{s10,s11,s12,s13,s14,s15}, complete synchronization \cite{s16}, generalized synchronization \cite{s17}, phase synchronization \cite{s18}, lag synchronization \cite{s19}, and so on.
\par Sometimes it takes a lot of people to achieve a goal, and different aims may be implemented by different people in real life, such as the division of research teams and the setting of government function departments. Because of different tasks, they have to achieve different goals. The same is true in neural networks. Each node in the same part wants to achieve the same goal, nodes in the different clusters want to achieve different goals, this is the so-called cluster synchronization. Cluster synchronization has important applications in science \cite{s20} and Bioengineering \cite{s21}. During the past few years, cluster synchronization has been widely studied by a large number of researchers, see \cite{s10,s11,s12,s13,s14,s15,s30,s31,s34,s36} and the relevant references. In an attempt to investigate cluster synchronization of complex systems, there have been many results related to cluster synchronization of coupled systems. This is not only because cluster synchronization has unusual correlations between network subsystems, but also because many coupled networks can exhibit excellent clustering characteristic. Additionally, As we all know, time delay is an ubiquitous phenomenon in physical complex networks. This is inevitable and may lead to a undesirable state of the system \cite{s35}. Therefore, the influence of time delay should be considered. As a result, it is essential to investigate cluster synchronization of coupled neural networks with time-varying delay.
\par Because of the complexity of the network structure, researchers have proposed many effective methods to synchronize the complex networks, some effective methods are pinning control \cite{s22,s23,s24}, impulsive control \cite{s25,s26,s27}, intermittent control \cite{s28}, adaptive control \cite{s29}, feedback control \cite{s36,s37,s38}, and so on. Impulse effect is a common phenomenon due to the instantaneous disturbance at some time. These phenomena are described by impulsive differential equation which have been used efficiently in modeling many practical problems that arise in the fields of engineering, physics, and science. Regarding to this, pinning impulsive control becomes a practically external control approach which is famous for low cost and energy-saving. Pinning impulsive control is characterized by adding impulsive controllers on only a fraction of network nodes. \cite{s26} investigated the pinning synchronization problem of stochastic impulsive networks using Lyapunov stability theory and pinning method. However, it only examined asymptotical synchronization and exponential synchronization of such dynamical networks in mean square. By introducing stochastic impulsive analysis technique and the selected nodes of the stochastic network are pinned at the impulse moment, authors in \cite{s30} discussed cluster synchronization of delayed stochastic neural networks under pinning impulsive control. A disadvantage of the above two papers is that they only considers the asymptotic behavior of the system as time approaches infinity.
\par In control theory, a considerable number of papers are focused on the utilization of pinning control to help complex networks achieve cluster synchronization. Based on the Lyapunov stability theory, the article \cite{s39} is devoted to study the cluster synchronization for network of a complex dynamical network of non-identical nodes with delayed and non-delayed coupling under a pinning control scheme. The finite-time and fixed-time cluster synchronization problem for complex networks with or without pinning control are studied in \cite{s40}. However, the systems in this paper are special coupled networks where use a concrete function $\rm sig^{p}(x)$. The major objective of \cite{s31} is to deal with cluster synchronization of coupled fractional-order neural networks with time delay, while the two control laws in this work is not pinning control but the control of all nodes in these networks. Based on it, we want to extend the idea of pinning control to control some of the key nodes in each cluster, and exert no control over the remaining nodes in the cluster, or only exert weak control over the remaining nodes, which is so-called hybrid control.
\par Based on above analysis, we study cluster synchronization of coupled neural networks with time-varying delay both on asymptotic cluster synchronization and on finite-time cluster synchronization. Unlike the literature \cite{s31} which is devoted to constant delay, we deal with time-varying delay, so the method in the literature \cite{s31} can't be used here. We utilize Lyapunov-Razumikhin method to deal with time-vary delay. Two novel sufficient conditions are obtained on asymptotic and finite-time cluster synchronization of the neural networks by exerting pinning impulsive control and hybrid control, respectively. The article is arranged as follows. In section 2, some preliminaries on asymptotic and finite-time cluster synchronization are provided. In section 3, we give some major results. In section 4, a numerical example is given to check the effectiveness of the consequences. The final section is the conclusion of this paper.

\section{Preliminaries}
\subsection{Model Description}
\par Consider the following system which consists of $N$ coupled neural networks with time-varying delay, the $i$th neural network is described by
\begin{equation}
\left\{\begin{array}{c}
\dot x_{i}(t)=-C_{i} x_{i}(t)+A_{i} f\left(x_{i}(t)\right)+B_{i} f\left(x_{i}(t-\tau(t))\right)+\sum \limits_{j=1}^{N} g_{i j} x_{j}(t)\\+I_{i}+u_{i}(t),\\
x_{i}(t)=\phi_{i}(t), t \in\left[-\tau, 0\right],
\end{array}\right.
\end{equation}
where $x_{i}(t)=\left(x_{i 1}(t), x_{i 2}(t), \cdots, x_{i n}(t)\right)^{T}$ is the state vector of the $i$th neural network at time $t$, $C_{i}>0$(means a symmetrical positive determined matrix) denotes the self-inhibition rate, $A_{i}=\left(a_{l r}^{i}\right)_{n \times n} \in R^{n \times n}$ and $B_{i}=\left(b_{l r}^{i}\right)_{n \times n} \in R^{n \times n}$ are connection weight matrices and delayed connection weight matrices, respectively. $f(x)$ is the activation function at time $t$ , $\tau(t)$ is transmittal delay. We suppose there exists two constants $\tau$ and $\sigma$, render $0<\tau(t)<\tau$ and $\dot {\tau}(t) \leq \sigma <1$. $I_{i} \in R^{n}$ means a input from the outside of the system, $u_{i}(t)$ is the controller that will be designed below. The element $g_{i j}$ of a linear coupling matrix is expressed as that if there is a connection from node $i$ to another node $j$ , then $g_{i j} \neq 0$, otherwise $g_{i j} =0$. The diagonal terms of the matrix satisfy the following diffusion coupling conditions: $g_{ii}=-\sum\limits_{j=1,j\ne i}^{N}g_{ji}$.

\subsection{Cluster Synchronization}
Let$\{\mathscr{C}_1,\cdots,\mathscr{C}_M\}$ be a partition of the set $\{1,\cdots,N\}$ , by which $\{1,\cdots,\\N\}$ is divided into $M$ nonempty subsets, i.e. $\mathscr{C}_l\ne \phi$ and $\bigcup \limits_{l=1}^{M}{\mathscr{C}_l}=\{1,\cdots,N\}$. For each $l\in\{1,\cdots,N\}$, let $\widetilde{l}$ denotes the subscript of the subset where the number $l$ is, i.e. $l \in \mathscr{C}_{\widetilde{l}}$.

\par A network with $N$ nodes is said to reach asymptotic cluster synchronization if the state variables of the network satisfies $\lim\limits _{t \rightarrow \infty}\left\|x_{i}(t)-x_{j}(t)\right\|=0$ for ${\widetilde{i} = \widetilde{j}}$ and $\lim \limits_{t \rightarrow \infty}\left\|x_{i}(t)-x_{j}(t)\right\|\ne 0$ for ${\widetilde{i} \ne \widetilde{j}}$. Furthermore, A network with $N$ nodes is said to reach finite-time cluster synchronization if there exists a $t^{*}>t_{0}$ such that for ${\widetilde{i} = \widetilde{j}}$ ,$\lim \limits_{t \rightarrow t^{*}}\left\|x_{i}(t)-x_{j}(t)\right\|=0$ and $\|x_{i}(t)-x_{j}(t)\|=0$ for $t \geq t^{*}$.
\\ \hspace*{\fill} \\
\textbf{Definition 1}\cite{s31}. The function $H=(h_{ij} \in R^{m\times m})$ is said to belong to class $\mathscr{A}_{1}$, if it satisfies $h_{i j} \geq 0,  h_{i i}=-\sum\limits_{j=1, j \neq i}^{m} h_{i j}$.
\\ \hspace*{\fill} \\
\textbf{Definition 2} \cite{s31}. The function $H=(h_{ij} \in R^{m\times m})$ belongs to class $\mathscr{A}_{2}$, if it satisfies $\sum\limits_{j=1}^{m}h_{i j}=0 $.
\\ \hspace*{\fill} \\
\textbf{Assumption 1} \cite{s31}. For $i \in \mathscr{C}_{p},p=\{1,2,\cdots,M\}$, there exists $\xi_{p}>0$, for all $x,y \in R^{n}$, $\|f (x)-f(y)\|\le \xi_{p} \|x-y\| $.
\\ \hspace*{\fill} \\
\textbf{Assumption 2} \cite{s31}. If all $N$ nodes are divided into $M$ clusters, The Laplace matrix can be defined in the following form
$$\begin{array}{l}
\left(\begin{array}{ccccc}
\bar{A}_{11} & \bar{A}_{12} & \bar{A}_{13} & \cdots & \bar{A}_{1 M} \\
\bar{A}_{21} & \bar{A}_{22} & \bar{A}_{23} & \cdots & \bar{A}_{2 M} \\
\bar{A}_{31} & \bar{A}_{32} & \bar{A}_{33} & \cdots & \bar{A}_{3 M} \\
\vdots & \vdots & \vdots & \ddots & \vdots \\
\bar{A}_{M 1} &\bar{A}_{M 2} & \bar{A}_{M 3} & \cdots & \bar{A}_{M M}
\end{array}\right),
\end{array}$$
where $\bar{A}_{i j} \in R^{\left(v_{i}-v_{i-1}\right) \times\left(v_{j}-v_{j-1}\right)},  \bar{A}_{p p} \in \mathscr{A}_{1}, q=1,2, \cdots, M $,
$ \bar{A}_{p q} \in \mathscr{A}_{2}, p, q=1,2, \cdots, M, p \neq q$.

\subsection{Mathematical Preparation}
\par Provided that the $N$ nodes can be divided into $M$ clusters, that is to say $\{1,2, \cdots, N\}=\mathscr{C}_{1} \cup \mathscr{C}_{2} \cup \cdots \cup \mathscr{C}_{M}$, where $\mathscr{C}_{p}=\left\{v_{p-1}+1, v_{p-1}+2, \cdots, v_{p}\right\}$, $p=1,2, \cdots, M$ and $v_{0}=0$, $v_{M}=N$, $v_{p-1}< v_{p}$. It is necessary to assume that the coefficients of system (1) meet $C_{v_{p-1}+1}=C_{v_{p-1}+2}=\cdots=C_{r_{p}}=C_{p}$,
$A_{p_{p-1}+1}=A_{p_{p-1}+2}=\cdots=A_{r_{p}}=A_{p}$, $B_{v_{p-1}+1}=B_{v_{p-1}+2}
=\cdots=B_{v_{p}}=B_{p}$, $I_{v_{p-1}+1}=I_{v_{p-1}+2}=\cdots=I_{v_{p}}=I_{p}$.
Suppose we want to control network (1) to achieve the following ideal state $x_{v_{0}+1}(t), \cdots, x_{v_{1}}(t) \rightarrow s_{1}(t)$;  $x_{v_{1}+1}(t), \cdots, x_{v_{2}}(t) \rightarrow s_{2}(t)$; $x_{v_{M}+1}(t), \cdots, x_{v_{M}}(t) \rightarrow s_{M}(t)$.
That is, $s=\left(s_{1}(t), \cdots, s_{1}(t), s_{2}(t), \cdots, s_{2}(t), \cdots, s_{M}(t), \cdots, s_{M}(t)\right) \subset R^{n \times N}$ is the desired leader.
$s_{p}(t)$can be defined as
\begin{equation}
\left. \dot s_{p}(t)=-C_{p} s_{p}(t)+A_{p} f\left(s_{p}(t)\right)+B_{p} f\left(s_{p}(t-\tau(t))\right)\right)+I_{p}, p=1,2, \cdots, M.
\end{equation}
Then there is the error dynamic system
\begin{equation}\begin{array}{*{20}{l}}
{\dot{e_i}(t) = \dot {x_i}(t) - \dot {s_p}(t)}\\
\begin{array}{l}
\;\;\;\;\;\; = -{C_p}{e_i}(t) + {A_p}\left[ {f\left( {{x_i}(t)} \right) - f\left( {{s_p}(t)} \right)} \right]\\
\;\;\;\;\;\;\;\;\;\; + {B_p}\left[ {f\left( {{x_i}(t - \tau (t))}\right) - f\left( {{s_p}(t - \tau (t))} \right)} \right]
+ \sum\limits_{j = 1}^N {{g_{ij}}} {e_j}(t) + {u_i}(t).
\end{array}
\end{array}
\end{equation}
\textbf{Definition 3} \cite{s32}. Let $V:R^{+}\times R^{n} \rightarrow R^{+}$, then $V$ is said to belong to class $\nu_{0}$, if \\
1.\begin{equation}
\lim\limits _{(t, y) \rightarrow\left(\tau_{k}^{+}, x\right)} V(t, y) = V\left(\tau_{k}^{+}, x\right),
\end{equation}
2. $V$ is locally Lipschitzian in $x$ and $V(t,0)=0$ for all $t \in R^{+}$.

Consider a general impulsive delayed differential equation(IDDE) in the form of
\begin{equation}\left\{\begin{array}{c}
\dot{x}(t)=f\left(t, x_{t}\right), t \neq t_{k}, t \geq t_{0}, \\
\Delta x\left(t_{k}\right)=I_{k}\left(t_{k}, x_{t_{k}}\right), k \in N, \\
x_{t_{0}}=\phi \in P C\left([-\tau, 0], R^{n}\right),
\end{array}\right.\end{equation}
where $x(t)=\left[x_{1}(t), x_{2}(t), \cdots, x_{n}(t)\right]^{T}$, $x_{t}=\{x(t+\theta): -\tau\le\theta\le0\}$. $f:R^{+}\times R^{n}\times P C\left([-\tau, 0], R^{n}\right)\rightarrow R^{n} $ is continuous. $\Delta x\left(t_{k}\right)=x(t_{k}^{+})-x(t_{k})$ denotes impulsive interval, $I:R^{+}\times R^{n} \rightarrow R^{n}$
is also continuous. We shall assume $f(t,0)=I_{k}(t,0)=0$ for all $t \in R^{+}$ and $k \in N $ so that
system (5) admits the trivial solution.
\\ \hspace*{\fill} \\
\textbf{Lemma 1} \cite{s32}. We suppose there exist $V\in\nu_{0}$ and $p>0$, $c>0$, $c_{1}>0$, $c_{2}>0$, $\sigma>0$, $\lambda >0$, $\gamma\ge1$ and $\sigma-\lambda\ge c$, make the following four conditions true\\
(i)$c_{1}\|x\|^{s} \leq V(t, x) \leq c_{2}\|x\|^{s}$ ,  $\forall t\in R^{+} $, $x \in R^{n}$;\\
(ii)$D^{+}V(t,\varphi (0)) \leq cV(t,\varphi (0))$, $\forall t\in [t_{k-1},t_{k}), k\in N$ hold, no matter when\\
$qV(t,\varphi (0))\ge V(t+\theta$, $\varphi (\theta))$ for $\theta\in [-\tau,0]$, where the constant
$q\ge \gamma e^{\lambda \tau}$;\\
(iii)$D^{+}V(t_{k}, \varphi (0)+I_{k}(t_{k},\varphi))\leq \eta_{k} V(t_{k}^{-}, \varphi (0))$, where
$0<\eta_{k-1}\leq1$, $k\in N$;\\
(iv)$\gamma \ge 1/\eta_{k-1}$ and $ln \eta_{k-1}<-(\sigma+\lambda)(t_{k}-t_{k-1}), k\in N$;\\
then IDDE(5) is global stability exponentially to the zero solution. The rate of convergence is $\frac{\lambda}{s}$.
\\ \hspace*{\fill} \\
\textbf{Lemma 2} \cite{s30}. $\forall x,y \in R_{n}$, there exist $\varpi>0$, $E \in R^{n \times n}>0$(positive definite metrix),
such that $2x^{T}y \le \varpi x^{T}Ex+\varpi^{-1}y^{T}E^{-1}y $.
\\ \hspace*{\fill} \\
\textbf{Lemma 3} \cite{s34}. Let $x_{1}, x_{2}, \cdots, x_{n} \in R^{n}$ are any vectors and $0<q<2$ is a real number satisfying:
$$\left(\left\|x_{1}\right\|^{2}+\left\|x_{2}\right\|^{2}+\cdots+\left\|x_{n}\right\|^{2}\right)^{q / 2}\le\left\|x_{1}\right\|^{q}+\left\|x_{2}\right\|^{q}+\cdots+\left\|x_{n}\right\|^{q} .$$
\\ \hspace*{\fill} \\
\textbf{Lemma 4} \cite{s34}. Suppose that a continuous, positive-definite function $V(t)$ satisfies $\dot V(t)\le-\alpha V^{\eta}(t), \forall t\ge t_{0}, V(t_{0})\ge 0,$ where $\alpha>0, 0<\eta<1$ are two constants, then for any given $t_{0}$, $V(t)$ satisfies $V^{1-\eta}(t)\le V^{1-\eta}(t_{0})-\alpha(1-\eta)(t-t_{0}), t_{0}\le t \le t_{1}$, and $V(t)=0, \forall t\ge t_{1}$, with $t_{1}$ given by $t_{1}=t_{0}+\frac{V^{1-\eta}(t_{0})}{\alpha(1-\eta)}$.

\section{Main Results}
\subsection{Asymptotic Cluster Synchronization}
In this part, we obtain a sufficient condition for asymptotic cluster synchronization of a class of coupled neural networks (1) with time-varying delays by means of pinning impulsive control.

\par In order to achieve the purpose of asymptotic cluster synchronization, pinning impulsive control only applies to certain nodes. We set $\delta (\cdot)$ as Dirac delta function, $d_{k}$ as impulse gain. If $-2<d_{k}<0$, we choose the first $\rho_{p}$ node to pin in the cluster $\mathscr{C}_{p}$, the selection method follows the form (a); if $d_{k} \ge 0$ or $d_{k}\le -2$, we choose the first $\rho_{p}$ node as the pinning node in the cluster $\mathscr{C}_{p}$, the selection method follows the form (b). We define the set of all pinned nodes as $\chi_{p}\left(t_{k}\right)=\left\{i_{p 1}, i_{p 2}, \cdots, i_{p \rho_{p}}\right\} \subset \mathscr{C}_{p}$. In order to add a pinning impulsive control, we set $t_{k} \geq 0$ are impulsive moments, and assume that $0=t_{0}<t_{1}< \cdots <t_{k}<t_{k+1}< \cdots$, $ \lim\limits_{k\to\infty}t_{k}=+\infty$, $\sup\limits _{k \geq 0}\left\{t_{k+1}-t_{k}\right\}<+\infty$. When $t=t_{k}$, the node errors in each cluster $\mathscr{C}_{p}$ are divided into the following two forms:\\
(a)$\|e_{i_{p 1}}(t_{k})\| \geq\|e_{i_{p 2}}(t_{k})\| \geq \cdots \geq\|e_{i_{p_{s}}}(t_{k})\| \geq\|e_{i_{p_{p(s+1)}}}(t_{k})\| \geq \cdots \geq\\ \|e_{i_{p_{k}(v_{p}-v_{p-1})}}(t_{k})\|;$\\
(b)$\|e_{i_{p 1}}(t_{k})\| \leq\|e_{i_{p 2}}(t_{k})\| \leq \cdots \leq\|e_{i_{p_{s s}}}(t_{k})\| \leq\|e_{i_{p_{p(s+1)}}(t_{k})}\| \leq \cdots \leq\\ \|e_{i_{p_{k}(v_{p}-p_{p-1})}}(t_{k})\|;$\\
where $i_{ps} \in \mathscr{C}_{p}, p=1,2,\cdots,M $, $i_{pu}\ne i_{pw}$, when $u \ne w$. What's more, providing $| e_{i_{p s}}\left(t_{k}\right)\|=\| e_{i_{p,(s+1)}}\left(t_{k}\right) \|$, then $i_{ps}<i_{p,s+1}$.

Design the following pinning impulsive controllers
\begin{equation}u_{i}(t)=\left\{\begin{array}{c}
\sum\limits_{k=1}^{+\infty} d_{k} e_{i}(t) \delta\left(t-t_{k}\right), i \in \chi_{p}\left(t_{k}\right), \\
0, i \in \mathscr{C}_{p} \backslash \chi_{p}\left(t_{k}\right).
\end{array}\right.\end{equation}
The form of the error dynamic system under the pinning impulsive controllers (6) is as follows
\begin{equation}\left\{\begin{array}{c}
\dot e_{i}(t)=-C_{p} e_{i}(t)+A_{p} F\left(e_{i}(t)\right)+B_{p} F\left(e_{i}(t-\tau(t))\right)+\sum\limits_{j=1}^{N} g_{i j} e_{j}(t),\\ \;\;\;\;\;\;\;\;\;\;\;\;\;\;\;\;\;\;\;\;\;\;\;\;\;\;\;\;\;\;\;\;\;\;\;\;\;\;\;\;\;\;\;\;\;\;\;\;\;\;i\in \mathscr{C}_{p},p=1,2,\cdots,M, t\ne t_{k},\\
e_{i}\left(t_{k}^{+}\right)=\left(1+d_{k}\right) e_{i}\left(t_{k}\right), i \in \chi_{p}\left(t_{k}\right), \\
e_{i}\left(t_{k}^{+}\right)=e_{i}\left(t_{k}\right), i \in \mathscr{C}_{p} \backslash \chi_{p}\left(t_{k}\right),
\end{array}\right.\end{equation}
where $F\left(e_{i}(t)\right)=f\left(x_{i}(t)\right)-f\left(s_{p}(t)\right)$, $F\left(e_{i}(t-\tau(t))\right)=f\left(x_{i}(t-
\tau(t))\right)-f\left(s_{p}(t-\tau(t))\right)$, $e_{i}\left(t_{k}^{+}\right)=\lim\limits _{h \rightarrow 0^{+}} e_{i}\left(t_{k}+h\right)$, $e_{i}\left(t_{k}^{-}\right)=\lim\limits _{h \rightarrow 0^{-}} e_{i}\left(t_{k}+h\right) $, $e_{i}(t)$ is left continuous when $t=t_{k}$. The initial condition of $e_{i}(t)$ is set as $e_{i}(t)=\phi(t) \in P C\left([-\tau, 0], R^{n}\right)$, where $P C\left([-\tau, 0], R^{n}\right)$ is the sum of all piecewise continuous functions $\phi$, $\phi$ has a supremum norm $\|\phi\|=\sup \limits_{-\tau \leq \theta \leq 0}\|\phi(\theta)\|$.
\\ \hspace*{\fill} \\
\textbf{Theorem 1.} According to the assumptions 1-2, if there have $E_{1} \in R^{n \times n}>0$, $E_{2} \in R^{n \times n}>0$, and $Q \in R^{n \times n}>0$, and $q\ge \gamma e^{\lambda \tau}$, $\alpha>0$, $\beta>0$, and $\mu,\upsilon$, so as to the following conditions are satisfied:
\begin{equation}\left(\begin{array}{ccc}
\gamma I_{N} & \sqrt{\alpha} I_{N} \otimes Q A_{p} & \sqrt{\beta} I_{N} \otimes Q B_{p} \\
* & -I_{N} \otimes E_{1} & 0 \\
* & * & -I_{N} \otimes E_{2}
\end{array}\right)<0,\end{equation}
\begin{equation}\begin{array}{l}
\beta^{-1} \xi^{2} E_{2}<v Q, \\
\end{array}\end{equation}
\begin{equation}\begin{array}{l}
\ln \eta_{k} \leq-(\sigma+\lambda) \Delta_{k-1}, \\
\end{array}\end{equation}
\begin{equation}\begin{array}{l}
\mu+\varepsilon \upsilon-(\sigma+\lambda)<0,
\end{array}\end{equation}
then the system (1) can realize cluster synchronization with exponential convergence rate $\frac{\lambda}{2}$ under the pinning impulsive controllers (6).
\\ \hspace*{\fill} \\
\textbf{Proof.} Set up a Lyapunov function as
\begin{equation}
V(t)=\sum\limits_{i=1}^{N} e_{i}^{T}(t) Q e_{i}(t),
\end{equation}
where $Q \in R^{n \times n}>0$.

It is easy to observe that $\lambda_{\min }(Q)\left\|e_{i}(t)\right\|^{2} \leq V(t)\leq \lambda_{\max }(Q)\left\|e_{i}(t)\right\|^{2},\\ t \in R^{+}$,
so condition (i) in Lemma 1 is satisfied and $s=2$.

For $t\in (t_{k-1},t_{k}]$, from the formulation (7), we can deduce
\begin{equation}\begin{aligned}
\frac{d V(t)}{d t}=2 & \sum_{p=1}^{M} \sum_{i=v_{p-1}}^{v_{p}} e_{i}^{T}(t) Q\left[-C_{p} e_{i}(t)+A_{p} F\left(e_{i}(t)\right)\right.\\
&+B_{p} F\left(e_{i}(t-\tau(t))\right)+\sum_{j=1}^{N} g_{i j} e_{j}(t)].
\end{aligned}\end{equation}

For $i\in \mathscr{C}_{p}$, according to the assumption 1 and lemma 2, exist $\alpha>0, \beta>0$ and $E_{1} \in R^{n \times n}>0$, $E_{2} \in R^{n \times n}>0$, such that
\begin{equation}\begin{array}{l}
2 e_{i}^{T}(t) Q A_{p} F\left(e_{i}(t)\right)
\leq \alpha e_{i}^{T}(t) Q A_{p} E_{1}^{-1} A_{p}^{T} Q e_{i}(t)+\alpha^{-1} F^{T}\left(e_{i}(t)\right) E_{1} F\left(e_{i}(t)\right) \\
\;\;\;\;\;\;\;\;\;\;\;\;\;\;\;\;\;\;\;\;\;\;\;\;\;\;\;\;\;\leq \alpha e_{i}^{T}(t) Q A_{p} E_{1}^{-1} A_{p}^{T} Q e_{i}(t)+\alpha^{-1} \xi_{p}^{2} e_{i}^{T}(t) E_{1} e_{i}(t),
\end{array}\end{equation}
together with
\begin{equation}\begin{array}{l}
2 e_{i}^{T}(t) Q B_{p} F\left(e_{i}(t-\tau(t))\right)  \leq \beta e_{i}^{T}(t) Q B_{p} E_{2}^{-1} B_{p}^{T} Q e_{i}(t)\\ +\beta^{-1} F^{T}\left(e_{i}(t-\tau(t))\right) \times E_{2} F\left(e_{i}(t-\tau(t))\right) \leq \beta e_{i}^{T}(t) Q B_{p} E_{2}^{-1} B_{p}^{T} Q e_{i}(t)\\+\beta^{-1} \xi_{p}^{2} e_{i}^{T}(t-\tau(t)) E_{2} e_{i}(t-\tau(t)).
\end{array}\end{equation}
Let $e(t)=(e_{1}^{T}(t),e_{2}^{T}(t),\cdots,e_{N}^{T}(t))^{T}$, then we have
\begin{equation}2 \sum\limits_{p=1}^{M} \sum\limits_{i=v_{p-1}+1}^{v_{p}} e_{i}^{T}(t) Q \sum\limits_{j=1}^{N} g_{i j} e_{j}(t)=2 e^{T}(t)(G \otimes Q)e(t).
\end{equation}
Substituting (14)-(16) into (13) and we can get
\begin{equation}\begin{array}{l}
\frac{d V(t)}{d t} =\sum\limits_{i=1}^{N} e_{i}^{T}(t) \Phi e_{i}(t)+\sum\limits_{p=1}^{M} \sum\limits_{i=v_{p-1}+1}^{v_{p}} e_{i}^{T}(t)\left(\alpha^{-1} \xi_{p}^{2} E_{1}\right) e_{i}(t)\\+2 e^{T}(t)(G \otimes Q) e(t)+\beta^{-1} \sum\limits_{p=1}^{M} \sum\limits_{i=v_{p-1}+1}^{v_{p}} e_{i}^{T}(t-\tau(t)) \xi_{p}^{2} Q e_{i}(t-\tau(t)),\\
\end{array}\end{equation}
where $\Phi=-Q C+C^{T} Q+\alpha Q A_{p} E_{1}^{-1} A_{p}^{T} Q+\beta Q B_{p} E_{2}^{-1} B_{p}^{T} Q$. According to (8) and the theory of schur complement, we can easily reason out
\begin{equation}\begin{array}{l}
\frac{d V(t)}{d t} \leq \mu \sum\limits_{i=1}^{N} e_{i}^{T}(t) Q e_{i}(t)+\upsilon[\sum\limits_{i=1}^{N} e_{i}^{T}(t-\tau(t)) Q e_{i}(t-\tau(t))] \\
\;\;\;\;\;\;\;\;\leq \mu V(t)+\upsilon[\sup\limits _{-\tau \leq \theta \leq 0} V(t+\theta)].\\
\end{array}\end{equation}

Due to $q\ge \gamma e^{\lambda \tau}$ is a constant, from Razumikhin condition, for all $t\in [t_{k-1},t_{k}),k\in N$, $qV(t,\varphi (0))\ge V(t+s,\varphi (s))$ for $s \in [-\tau,0]$, we have $e_{i}^{T}(t-\tau(t)) Q e_{i}(t-\tau(t)) \le qe_{i}^{T}(t)Q e_{i}(t), t\in [t_{k-1},t_{k}),k\in N$. So the formula (18) can be further deduced as
\begin{equation}\begin{array}{l}
\frac{d V(t)}{d t} \leq (\mu + \upsilon q)V(t).
\end{array}\end{equation}
Therefore, the parameter $c$ in condition (ii) in Lemma 1 can be chosen as $\mu + \upsilon q$, so that condition (ii) can be realized.

Next, let us consider when $t=t_{k}$, i.e. $e_{i}\left(t_{k}^{+}\right)=\left(1+d_{k}\right) e_{i}\left(t_{k}\right), i \in \chi_{p}\left(t_{k}\right); e_{i}\left(t_{k}^{+}\right)=e_{i}\left(t_{k}\right), i \in \mathscr{C}_{p} \backslash \chi_{p}\left(t_{k}\right)$.
We can find that
\begin{equation}\begin{array}{l}
V\left(t_{k}^{+}\right)=\sum\limits_{i=1}^{N} e_{i}^{T}\left(t_{k}^{+}\right) Q e_{i}\left(t_{k}^{+}\right) \\
\;\;\;\;\;\;\;\;\;\;\;=\sum\limits_{p=1}^{M} \sum\limits_{i \in \chi_{p}\left(t_{k}\right)} e_{i}^{T}\left(t_{k}^{+}\right) Q e_{i}\left(t_{k}^{+}\right)+\sum\limits_{p=1}^{M} \sum\limits_{i \in \mathscr{C}_{p} \backslash \chi_{p}\left(t_{k}\right)} e_{i}^{T}\left(t_{k}^{+}\right) Q e_{i}\left(t_{k}^{+}\right) \\
\;\;\;\;\;\;\;\;\;\;\;=\left(1+d_{k}\right)^{2} \sum\limits_{p=1}^{M} \sum\limits_{i \in \chi_{p}\left(t_{k}\right)} e_{i}^{T}\left(t_{k}\right) Q e_{i}\left(t_{k}\right)-d_{k}\left(2+d_{k}\right) \\
\;\;\;\;\;\;\;\;\;\;\;\;\;\;\times \sum\limits_{p=1}^{M} \sum\limits_{i \in \mathscr{C}_{p} \backslash \chi_{p}\left(t_{k}\right)} e_{i}^{T}\left(t_{k}\right) Q e_{i}\left(t_{k}\right).
\end{array}\end{equation}
If $-2 < d_{k} < 0$, in the light of pinned nodes in $\chi_{p}\left(t_{k}\right)$, we have
\begin{equation}\begin{array}{l}
\frac{1}{v_{p}-v_{p-1}-\rho_{p}} \sum\limits_{i \in \mathscr{C}_{p} \backslash \chi_{p}\left(t_{k}\right)} e_{i}^{T}\left(t_{k}\right) Q e_{i}\left(t_{k}\right) \\
\;\;\;\;\;\;\;\;\;\;\;\;\leq \frac{\lambda_{\max }(Q)}{v_{p}-v_{p-1}-\rho_{p}} \sum\limits_{i \in \mathscr{C}_{p} \backslash \chi_{p}\left(t_{k}\right)} e_{i}^{T}\left(t_{k}\right) e_{i}\left(t_{k}\right) \\
\;\;\;\;\;\;\;\;\;\;\;\;\leq \frac{\lambda_{\max }(Q)}{v_{p}-v_{p-1}} \sum\limits_{i=v_{p-1}+1}^{v_{p}} e_{i}^{T}\left(t_{k}\right) e_{i}\left(t_{k}\right) \\
\;\;\;\;\;\;\;\;\;\;\;\;\leq \frac{\lambda_{\max }(Q)}{\lambda_{\min }(Q)\left(v_{p}-v_{p-1}\right)} \sum\limits_{i=v_{p-1}+1}^{v_{p}} e_{i}^{T}\left(t_{k}\right) Q e_{i}\left(t_{k}\right).
\end{array}\end{equation}
So we can get
\begin{equation}\begin{array}{l}
V\left(t_{k}^{+}\right) \leq\left(1+d_{k}\right)^{2} \sum\limits_{p=1}^{M} \sum\limits_{i=v_{p-1}+1}^{v_{p}} e_{i}^{T}\left(t_{k}\right) Q e_{i}\left(t_{k}\right)-d_{k}\left(2+d_{k}\right) \\
\;\;\;\;\;\;\;\;\;\;\;\;\;\;\;\;\;\;\;\times \sum\limits_{p=1}^{M} \frac{\lambda_{\max }(Q)\left(v_{p}-v_{p-1}-\rho_{p}\right)}{\lambda_{\min }(Q)\left(v_{p}-v_{p-1}\right)} \sum\limits_{i=v_{p-1}+1}^{v_{p}} e_{i}^{T}\left(t_{k}\right) Q e_{i}\left(t_{k}\right) \\
\;\;\;\;\;\;\;\;\;\;\;\;= \left[\left(1+d_{k}\right)^{2}-d_{k}\left(2+d_{k}\right) \frac{\lambda_{\max }(Q)\left(v_{p}-v_{p-1}-\rho_{p}\right)}{\lambda_{\min }(Q)\left(v_{p}-v_{p-1}\right)}\right] \sum\limits_{i=1}^{N} e_{i}^{T}\left(t_{k}\right) Q e_{i}\left(t_{k}\right) \\
\;\;\;\;\;\;\;\;\;\;\;\; \doteq \eta_{k} V\left(t_{k}\right).
\end{array}\end{equation}

The same can be found if $d_{k}\le -2$ or $d_{k}\ge 0$, so we omit here. Utilizing Lemma 1 and (9)-(11), the trivial solution of system (7) can realize globally exponentially stable with convergence rate $\frac{\lambda}{2}$. Hence, the neural network (1) achieves exponential cluster synchronization asymptotically with exponent $\frac\lambda{2}$ under the pinning impulsive controllers (6).

\subsection{Finite-time Cluster Synchronization}
\par Motivated by pinning control scheme, in this part, we obtain a sufficient condition for finite-time cluster synchronization of a class of coupled neural networks (1) with time-varying delays by means of hybrid control.

\par In order to achieve the goal of finite-time cluster synchronization, we redefine the notation $\chi_{p}$ as the nodes in the $p$th cluster which is directly connected with the nodes in other clusters. So the notation $\mathscr{C}_{p}\backslash \chi_{p}$ means the remaining nodes in the $p$th cluster which isn't directly connected to nodes in other clusters.

\par Contrive the following hybrid controllers
\begin{equation}
u_{i}(t)=\left\{\begin{array}{c}
-g_{1} e_{i}(t)-\phi_{i}(t), i \in \chi_{p},\\
-\varphi_{i}(t), i \in \mathscr{C}_{p} \backslash \chi_{p},
\end{array}\right.\end{equation}
where $\phi_{i}(t)=k{\rm sign}\left(e_{i}(t)\right)\left|e_{i}(t)\right|^{\mu}+k\left(k_{1} \int_{t-\tau(t)}^{t} e_{i}^{T}(s) e_{i}(s) d s\right)^{\frac{1+\mu}{2}} \times \\ \psi\left(e_{i}(t),\|e(t)\|\right)$, $\varphi_{i}(t)=k{\rm sign}\left(e_{i}(t)\right)\left|e_{i}(t)\right|^{\mu}+k\left(k_{1} \int_{t-\tau(t)}^{t} e_{i}^{T}(s) e_{i}(s) d s\right)^{\frac{1+\mu}{2}} \times \\\psi\left(e_{i}(t),\|e(t)\|\right)+2 k\left(g_{1} \int_{t}^{t_{1}} e_{i}^{T}(s) e_{i}(s) d s\right)^{\frac{1+\mu}{2}} \psi\left(e_{i}(t),\|e(t)\|\right)$, $i=1,2,\cdots,N$, $k$ is tunable, $g_{1}$ is a positive constant which has the practical meaning of control strength, and $k_{1}>0$. $\psi\left(e_{i}(t),\|e(t)\|\right)$ is denoted as follows: if $\|e(t)\| \ne0,  \psi\left(e_{i}(t),\|e(t)\|\right)=\frac{e_{i}(t)}{\|e(t)\|^{2}}$, or else $\psi\left(e_{i}(t),\|e(t)\|\right)=0$.
\par The form of the error dynamic system under the hybrid controllers (23) is as follows
\begin{equation}\left\{\begin{array}{c}
\dot{e}_{i}(t)=-C_{p} e_{i}(t)+A_{p} F\left(e_{i}(t)\right)+B_{p} F\left(e_{i}(t-\tau(t))\right)+\sum\limits_{j=1}^{N} g_{i j} e_{j}(t)\\-g_{1} e_{i}(t)-\phi_{i}(t),i \in \chi_{p},\\
\dot{e}_{i}(t)=-C_{p} e_{i}(t)+A_{p} F\left(e_{i}(t)\right)+B_{p} F\left(e_{i}(t-\tau(t))\right)+\sum\limits_{j=1}^{N} g_{i j} e_{j}(t)\\-\varphi_{i}(t),i \in \mathscr{C}_{p} \backslash \chi_{p}.
\end{array}\right.\end{equation}

\par Our next goal is to drive the coupled neural network (1) to achieve synchronization among the node states $x_{i}(t)$ and the objective state $s_{p}(t)$ over a finite-time interval by designing the hybrid controllers (23) and applying the finite-time stability theory.
\\ \hspace*{\fill} \\
\textbf{Theorem 2.} Suppose that Assumptions 1-2 hold, if the following conditions are satisfied:\\
(1)$\beta^{-1} \xi_{p}^{2} \lambda_{\max }\left(E_{2}\right)+\sigma<1$,\\
(2)$-\lambda_{\min }\left(C_{p}\right)+\frac{1}{2} \alpha \lambda_{\max }\left(A_{p} E_{1}^{-1} A_{p}^{T}\right)+\frac{1}{2} \alpha^{-1} \xi_{p}^{2} \lambda_{\max }\left(E_{1}\right)+\frac{1}{2} \beta \lambda_{\max }\left(B_{p} E_{2}^{-1} B_{p}^{T}\right)\\+\lambda_{\max }(G \otimes I)+\frac{1}{2} k_{1}-g_{1}<0$,\\
then by the hybrid controllers (23), the system (1) can achieve cluster synchronization in finite time and the settling time is $t_{1}=t_{0}+\frac{1}{k(1-\mu)}(2V(t_{0}))^{\frac{1-\mu}{2}}$.
\\ \hspace*{\fill} \\
\textbf{proof.} Set up a Lyapunov functional as
\begin{equation}\begin{aligned}
V(t)=\frac{1}{2} & \sum_{i=1}^{N} e_{i}^{T}(t) e_{i}(t)+\frac{1}{2} k_{1} \sum_{i=1}^{N} \int_{t-\tau(t)}^{t} e_{i}^{T}(s) e_{i}(s) d s \\
&+g_{1} \sum_{p=1}^{M} \sum_{\mathscr{C}_{p} \backslash \chi_{p}} \int_{t}^{t_{1}} e_{i}^{T}(s) e_{i}(s) d s.
\end{aligned}\end{equation}

Calculating the derivative of $V(t)$ along the solutions of system (24) outputs the following:
\begin{equation}\begin{aligned}
\dot{V}(t)=& \sum_{i=1}^{N} e_{i}^{T}(t)\left[-C_{p} e_{i}(t)+A_{p} F\left(e_{i}(t)\right)+B_{p} F\left(e_{i}(t-\tau(t))\right)\right.\\
&+\sum_{j=1}^{N} g_{i j} e_{j}(t)+u_{i}(t)]+\frac{1}{2} k_{1} \sum_{i=1}^{N} e_{i}^{T}(t) e_{i}(t) \\
&-\frac{1-i(t)}{2} k_{1} \sum_{i=1}^{N} e_{i}^{T}(t-\tau(t)) e_{i}(t-\tau(t)) \\
&-g_{1} \sum_{p=1}^{M} \sum_{i \in \mathscr{C}_{p} \backslash \chi_{p}} e_{i}^{T}(t) e_{i}(t).
\end{aligned}\end{equation}

Inserting the hybrid controllers (23) into $\dot V(t)$, we have
$$\begin{aligned}
\dot{V}(t)=& \sum_{i=1}^{N} e_{i}^{T}(t)\left[-C_{p} e_{i}(t)+A_{p} F\left(e_{i}(t)\right)+B_{p} F\left(e_{i}(t-\tau(t))\right)\right.\\
&+\sum_{j=1}^{N} g_{i j} e_{j}(t)]+\frac{1}{2} k_{1} \sum_{i=1}^{N} e_{i}^{T}(t) e_{i}(t)-\frac{1-\dot{\tau}(t)}{2} \times \\
& k_{1} \sum_{i=1}^{N} e_{i}^{T}(t-\tau(t)) e_{i}(t-\tau(t))-g_{1} \sum_{p=1}^{M} \sum_{i \in \mathscr{C}_{p} \backslash \chi_{p}} e_{i}^{T}(t) e_{i}(t) \\
&-\sum_{p=1}^{M} \sum_{i \in \chi_{p}} g_{1} e_{i}^{T}(t) e_{i}(t)-k \sum_{i=1}^{N} e_{i}^{T}(t)\left[\operatorname{sign}\left(e_{i}(t)\right)\left|e_{i}(t)\right|^{\mu}\right.\\
&+(k_{1} \int_{t-\tau(t)}^{t} e_{i}^{T}(s) e_{i}(s) d s)^{\frac{1+\mu}{2}} \psi\left(e_{i}(t),\|e(t)\|\right)] \\
&-2 k \sum_{p=1}^{M} \sum_{i \in \mathscr{C}_{p} \backslash \chi_{p}} e_{i}^{T}(t)[(g_{1} \int_{t}^{t_{1}} e_{i}^{T}(s) e_{i}(s) d s)^{\frac{1+\mu}{2}} \times\\
&\left.\psi\left(e_{i}(t),\|e(t)\|\right)\right],
\end{aligned}$$
that is
$$\begin{aligned}
\dot{V}(t)=&-\sum_{i=1}^{N} e_{i}^{T}(t) C_{p} e_{i}(t)+\sum_{i=1}^{N} e_{i}^{T}(t) A_{p} F\left(e_{i}(t)\right) \\
&+\sum_{i=1}^{N} e_{i}^{T}(t) B_{p} F\left(e_{i}(t-\tau(t))\right)+\sum_{i=1}^{N} e_{i}^{T}(t) \sum_{j=1}^{N} g_{i j} e_{j}(t) \\
&+\frac{1}{2} k_{1} \sum_{i=1}^{N} e_{i}^{T}(t) e_{i}(t)-\frac{1-\dot{\tau}(t)}{2} \times \\
& k_{1} \sum_{i=1}^{N} e_{i}^{T}(t-\tau(t)) e_{i}(t-\tau(t))-g_{1} \sum_{p=1}^{M} \sum_{i \in \mathscr{C}_{p} \backslash \chi_{p}} e_{i}^{T}(t) e_{i}(t) \\
&-\sum_{p=1}^{M} \sum_{i \in \chi_{p}} g_{1} e_{i}^{T}(t) e_{i}(t)-k \sum_{i=1}^{N} e_{i}^{T}(t)\left[\operatorname{sign}\left(e_{i}(t)\right)\left|e_{i}(t)\right|^{\mu}\right.\\
&+(k_{1} \int_{t-\tau(t)}^{t} e_{i}^{T}(s) e_{i}(s) d s)^{\frac{1+\mu}{2}} \psi\left(e_{i}(t),\|e(t)\|\right)] \\
&-2 k \sum_{p=1}^{M} \sum_{i \in \mathscr{C}_{p} \backslash \chi_{p}} e_{i}^{T}(t)[(g_{1} \int_{t}^{t_{1}} e_{i}^{T}(s) e_{i}(s) d s)^{\frac{1+\mu}{2}} \times\\
&\psi\left(e_{i}(t),\|e(t)\|\right)].
\end{aligned}$$
Utilizing the same technique as in the proof of Theorem 1 and $\dot\tau(t)\le \sigma$,
$$\begin{aligned}
\dot{V}(t) \leq &-\lambda_{\min }\left(C_{p}\right) e^{T}(t) e(t)+\frac{1}{2} \sum_{i=1}^{N}[\alpha e_{i}^{T}(t) A_{p} E_{1}^{-1} A_{p}^{T} e_{i}(t)\\
&+\alpha^{-1} \xi_{p}^{2} e_{i}^{T}(t) E_{1} e_{i}(t)]+\frac{1}{2} \sum_{i=1}^{N}[\beta e_{i}^{T}(t) B_{p} E_{2}^{-1} B_{p}^{T} e_{i}(t)\\
&+\beta^{-1} \xi_{p}^{2} e_{i}^{T}(t-\tau(t)) E_{2} e_{i}(t-\tau(t))]+e^{T}(t)(G \otimes I) e(t) \\
&+\frac{1}{2} k_{1} \sum_{i=1}^{N} e_{i}^{T}(t) e_{i}(t)-\frac{1-\sigma}{2} k_{1} \sum_{i=1}^{N} e_{i}^{T}(t-\tau(t)) e_{i}(t-\tau(t)) \\
&-g_{1} e^{T}(t) e(t)-k \sum_{i=1}^{N} e_{i}^{T}(t)[\operatorname{sign}(e_{i}(t))\left|e_{i}(t)\right|^{\mu}\\
&+(k_{1} \int_{t-\tau(t)}^{t} e_{i}^{T}(s) e_{i}(s) d s)^{\frac{1+\mu}{2}} \psi(e_{i}(t),\|e(t)\|)] \\
&-2 k \sum_{p=1}^{M} \sum_{i \in \mathscr{C}_{p} \backslash \chi_{p}} e_{i}^{T}(t)[(g_{1} \int_{t}^{t_{1}} e_{i}^{T}(s) e_{i}(s) d s)^{\frac{1+\mu}{2}} \times\\
&\psi(e_{i}(t),\|e(t)\|)].
\end{aligned}$$
Expanding this out yields
$$\begin{aligned}
\dot{V}(t) \leq &-\lambda_{\min }\left(C_{p}\right) e^{T}(t) e(t)+\frac{1}{2} \alpha \lambda_{\max }\left(A_{p} E_{1}^{-1} A_{p}^{T}\right) e^{T}(t) e(t) \\
&+\frac{1}{2} \alpha^{-1} \xi_{p}^{2} \lambda_{\max }\left(E_{1}\right) e^{T}(t) e(t)+\frac{1}{2} \beta \lambda_{\max }\left(B_{p} E_{2}^{-1} B_{p}^{T}\right) e^{T}(t) e(t) \\
&+\frac{1}{2} \beta^{-1} \xi_{p}^{2} \lambda_{\max }\left(E_{2}\right) e^{T}(t-\tau(t)) e(t-\tau(t)) \\
&+\lambda_{\max }(G \otimes I) e^{T}(t) e(t)+\frac{1}{2} k_{1} e^{T}(t) e(t) \\
&-\frac{1-\sigma}{2} k_{1} e^{T}(t-\tau(t)) e(t-\tau(t))-g_{1} e^{T}(t) e(t) \\
&-k \sum_{i=1}^{N} e_{i}^{T}(t)[\operatorname{sign}(e_{i}(t))|e_{i}(t)|^{\mu}+(k_{1} \int_{t-\tau(t)}^{t} e_{i}^{T}(s) e_{i}(s) d s)^{\frac{1+\mu}{2}}\\
&\left.\times \psi\left(e_{i}(t),\|e(t)\|\right)\right]-2 k \sum_{p=1}^{M} \sum_{i \in \mathscr{C}_{p} \backslash \chi_{p}} e_{i}^{T}(t)[(g_{1} \int_{t}^{t_{1}} e_{i}^{T}(s) e_{i}(s) d s)^{\frac{1+\mu}{2}}\\
& \times \psi(e_{i}(t),\|e(t)\|)].
\end{aligned}$$

By sorting through it, we have
$$\begin{aligned}
\dot{V}(t) \leq &[-\lambda_{\min }\left(C_{p}\right)+\frac{1}{2} \alpha \lambda_{\max }\left(A_{p} E_{1}^{-1} A_{p}^{T}\right)+\frac{1}{2} \alpha^{-1} \xi_{p}^{2} \lambda_{\max }\left(E_{1}\right)\\
&+\frac{1}{2} \beta \lambda_{\max }\left(B_{p} E_{2}^{-1} B_{p}^{T}\right)+\lambda_{\max }(G \otimes I)+\frac{1}{2} k_{1}-g_{1}] e^{T}(t) e(t) \\
&+(\frac{1}{2} \beta^{-1} \xi_{p}^{2} \lambda_{\max }\left(E_{2}\right)-\frac{1-\sigma}{2}) e^{T}(t-\tau(t)) e(t-\tau(t)) \\
&-k[\sum_{i=1}^{N} e_{i}^{T}(t)\left(\operatorname{sign}(e_{i}(t)\right)\left|e_{i}(t)\right|^{\mu}+k_{1} \int_{t-\tau(t)}^{t} e_{i}^{T}(s) e_{i}(s) ds)^{\frac{1+\mu}{2}}\\
&\times \psi(e_{i}(t),\|e(t)\|)]-2 k \sum_{p=1}^{M} \sum_{i \in \mathscr{C}_{p}} e_{i}^{T}(t)[(g_{1} \int_{t}^{t_{1}} e_{i}^{T}(s) e_{i}(s) d s)^{\frac{1+\mu}{2}}\\
&\times \psi(e_{i}(t),\|e(t)\|)].
\end{aligned}$$

From the conditions (1) and (2) of Theorem 2, we get that
$$\begin{aligned}
\dot{V}(t) \leq &-k \sum_{i=1}^{N} e_{i}^{T}(t)[\operatorname{sign}\left(e_{i}(t)\right)\left|e_{i}(t)\right|^{\mu}+(k_{1} \int_{t-\tau(t)}^{t} e_{i}^{T}(s) e_{i}(s) d s)^{\frac{1+\mu}{2}}\\
&\left.\psi\left(e_{i}(t),\|e(t)\|\right)\right]-2 k \sum_{p=1}^{M} \sum_{i \in \mathscr{C}_{p} \backslash \chi_{p}} e_{i}^{T}(t)[(g_{1} \int_{t}^{t_{1}} e_{i}^{T}(s) e_{i}(s) d s)^{\frac{1+\mu}{2}}\\
&\left.\times \psi\left(e_{i}(t),\|e(t)\|\right)\right] \\
\leq-& k \sum_{i=1}^{N}\left|e_{i}^{T}(t) e_{i}(t)\right|^{1+\mu}-k \sum_{i=1}^{N}(k_{1} \int_{t-\tau(t)}^{t} e_{i}^{T}(s) e_{i}(s) d s)^{\frac{1+\mu}{2}} \\
&-2 k \sum_{p=1}^{M} \sum_{\mathscr{C}_{p} \backslash \chi_{p}}(g_{1} \int_{t}^{t_{1}} e_{i}^{T}(s) e_{i}(s) d s)^{\frac{1+\mu}{2}}.
\end{aligned}$$

Through Lemma 3, one gain

$$\begin{aligned}
\dot{V}(t)\leq-& 2^{\frac{1+\mu}{2}} k(\frac{1}{2} \sum_{i=1}^{N} e_{i}^{2}(t)+\frac{1}{2} \sum_{i=1}^{N} k_{1} \int_{t-\tau(t)}^{t} e_{i}^{T}(s) e_{i}(s) d s\\
&+\sum_{p=1}^{M} \sum_{i \in \mathscr{C}_{p} \backslash \chi_{p}} g_{1} \int_{t}^{t_{1}} e_{i}^{T}(s) e_{i}(s) d s)^{\frac{1+\mu}{2}} \\
=&-2^{\frac{1+\mu}{2}} k V^{\frac{1+\mu}{2}}(t).
\end{aligned}$$

Utilizing Lemma 4, the error system (24) can achieve finite-time synchronization no matter what initial condition is, the settling time is $$t_{1}=t_{0}+\frac{1}{k(1-\mu)}(2V(t_{0}))^{\frac{1-\mu}{2}}.$$ Therefore, the neural network (1) can achieve finite-time cluster synchronization under the hybrid controllers (23).

\section{Illustrative example}
In this section, an example with three cases is presented to demonstrate the effectiveness of our results.

\par Example: Consider two-neuron time-varying delayed coupled neural networks (2), where $f(r)=(f_{1}(r),f_{2}(r))^{T}=(\rm arctan(r), \rm arctan(r))^{T}$, $\tau(t)=\frac{1.7e^{t}}{1+e^{t}}$, $M=2$, and $C_{1}=C_{2}=\left(\begin{array}{ll}1 & 0 \\ 0 & 1\end{array}\right)$, $A_{1}=\left(\begin{array}{cc}1.95 & -0.1 \\ -5 & 3\end{array}\right)$, $I_{1}=I_{2}=\left(\begin{array}{l}0 \\ 0 \end{array}\right)$, \\$A_{2}=\left(\begin{array}{cc}2 & -0.11 \\ -5.1 & 3\end{array}\right)$, $B_{1}=\left(\begin{array}{cc}-1.5 & -0.1 \\ -0.3 & -2.41\end{array}\right)$, $B_{2}=\left(\begin{array}{cc}-1.5 & -0.1 \\ -0.2 & -2.45\end{array}\right)$. The initial conditions are $S_{p}(t)=(0.4,0.6)^{T}$, $t\in [-1,0]$ with $p=1,2$.

\par Given a seven coupled time-varying delayed neural networks (1), the coupling matrix of which is
$$G=\left(\begin{array}{ccccccc}
-4 & 2 & 2 & 0.5 & -0.5 \\
2 & -3 & 1 & 0.2 & -0.2 \\
1 & 1 & -2 & 0.3 & -0.3 \\
0.5 & 0.5 & -1 & -1 & 1 \\
0.4 & -0.4 & 0 & 1 & -1 \\
\end{array}\right).$$

Then the system (1) can be divided into two clusters $\mathscr{C}_{1}=\{1,2\}, \mathscr{C}_{2}=\{3,4,5\}$.
\par Case 1: We set the error $\left\|e_{1}\right\|$, $\left\|e_{2}\right\|$ of synchronization with initial value $(10,-5),i=1,2$, $(8,-6),i=3,4,5$, respectively. The error dynamic system with $u_{i}(t)=0$ can't realize synchronization, which is depicted in Fig.1 and Fig.2.

\par Case 2: The error dynamic system under the pinning impulsive controllers (6) with $p_{1}=1, p_{2}=3, t_{k}=0.03k, d_{k}=-0.8$ can realize asymptotic cluster synchronization, which is depicted in Fig.3 and Fig.4.

\par Case 3: The error dynamic system under the hybrid controllers (23) with$k=2, k_{1}=1.4, g_{1}=41.4$ can realize finite-time cluster synchronization, which is depicted in Fig.5 and Fig.6.

\begin{figure}[htbp]
\begin{minipage}[t]{0.45\linewidth}
\centering
\includegraphics[height=5cm,width=5cm]{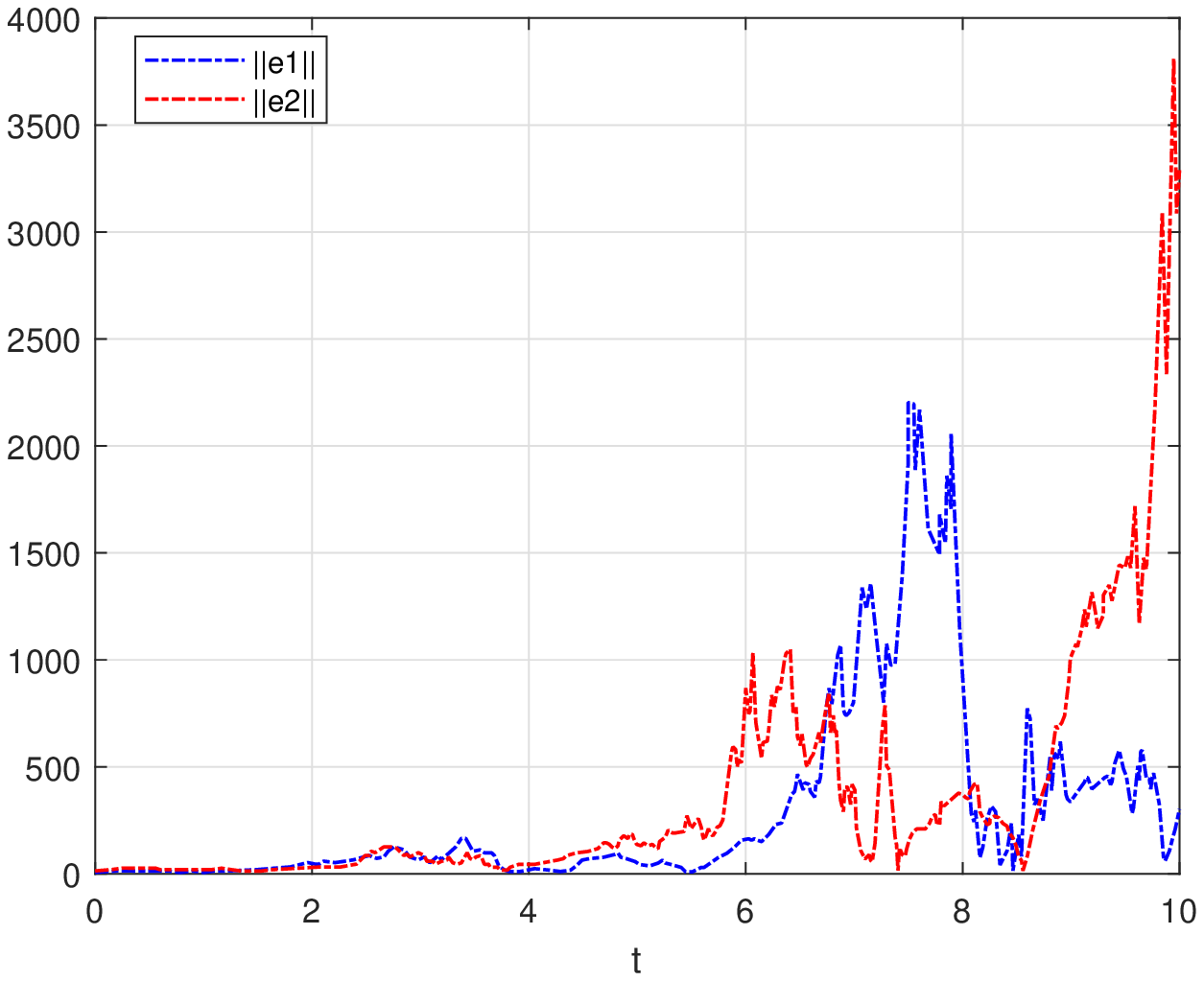}
\caption{The error dynamic trajectory when $i=1,2$ without controller.}
\end{minipage}
\begin{minipage}[t]{0.45\linewidth}
\centering
\includegraphics[height=5cm,width=5cm]{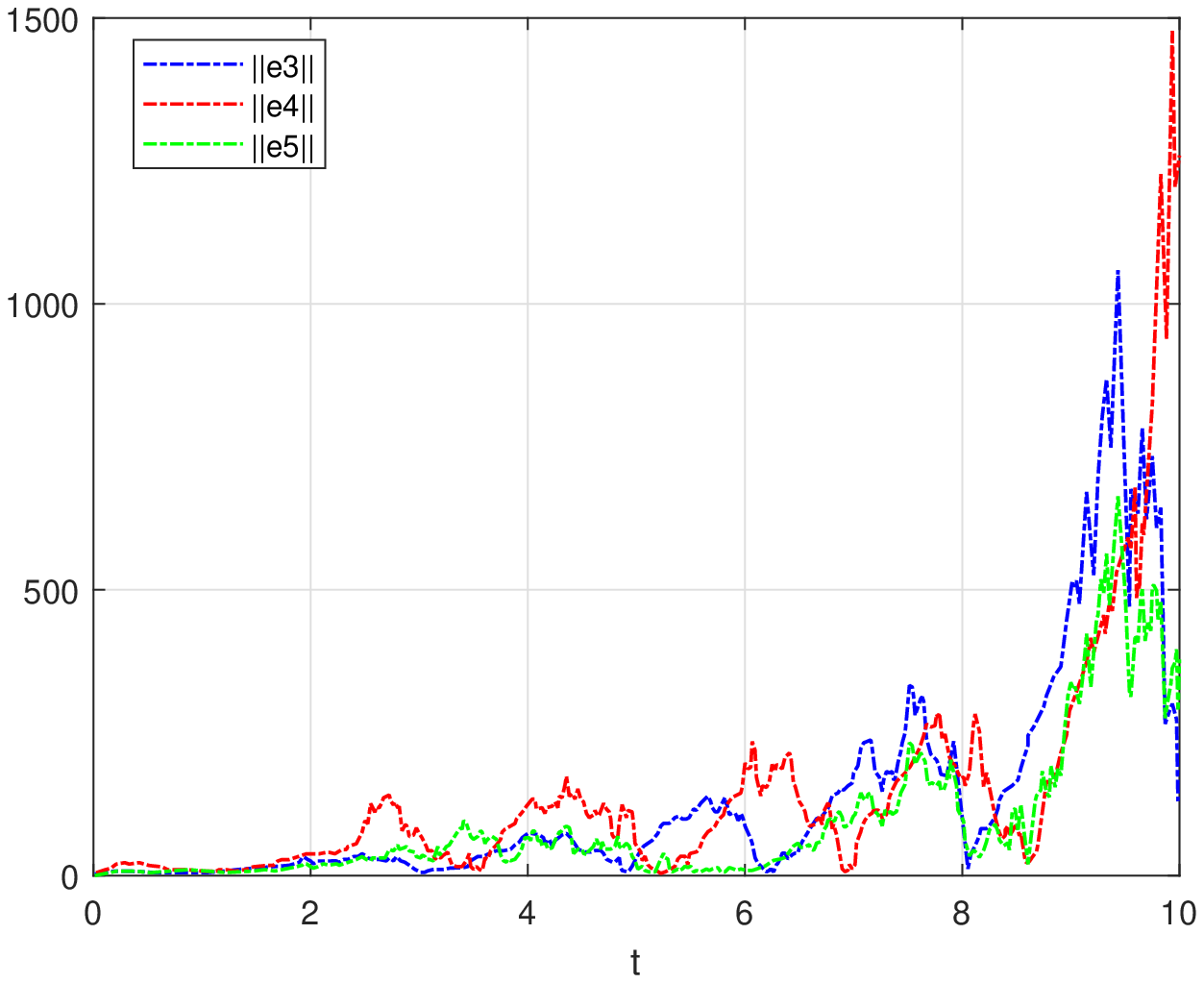}
\caption{The error dynamic trajectory when $i=3,4,5$ without controller.}
\end{minipage}
\end{figure}

\begin{figure}[htbp]
\begin{minipage}[t]{0.45\linewidth}
\centering
\includegraphics[height=5cm,width=5cm]{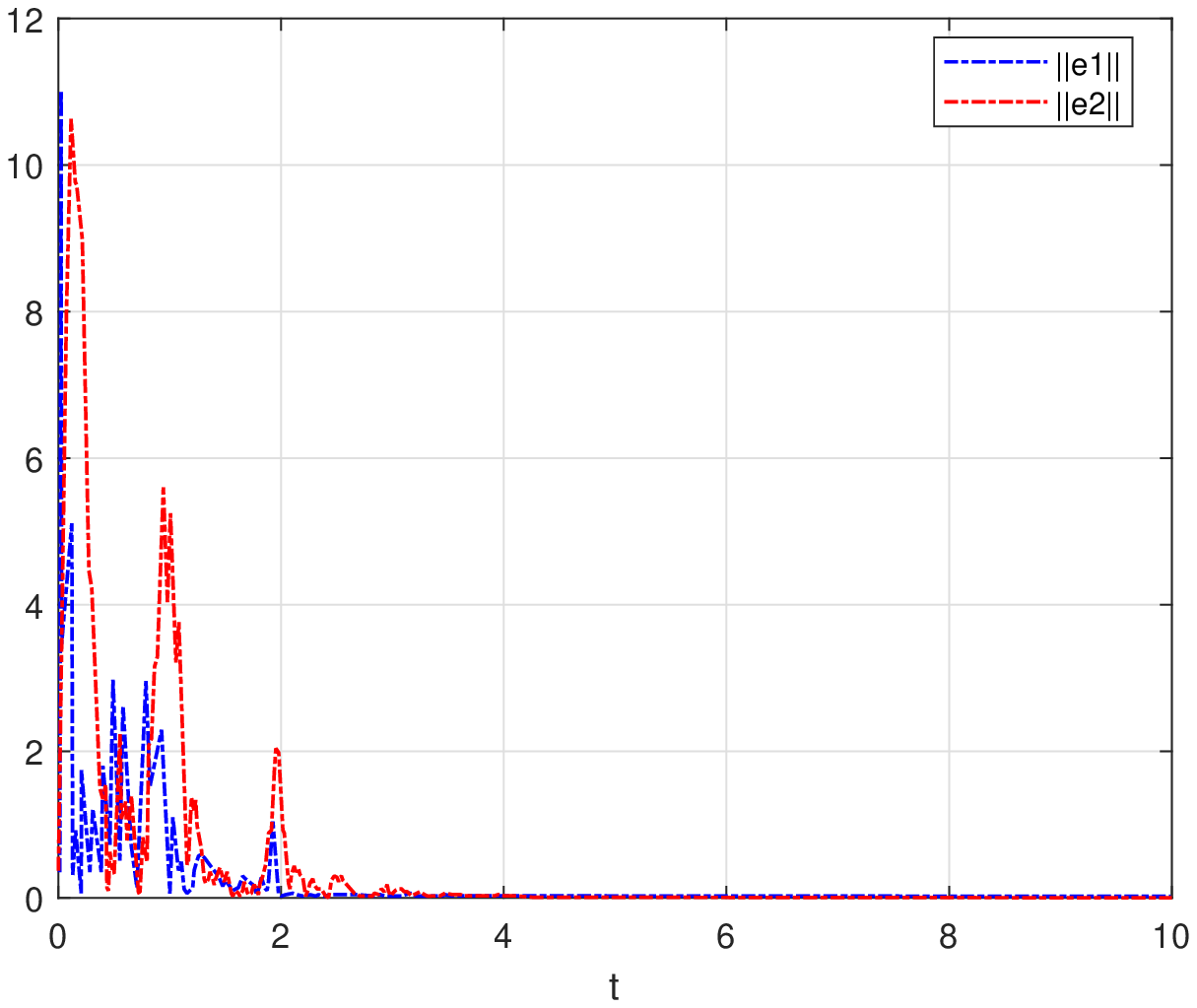}
\caption{The error dynamic trajectory when $i=1,2$ with the pinning impulsive controllers (6).}
\end{minipage}
\begin{minipage}[t]{0.45\linewidth}
\centering
\includegraphics[height=5cm,width=5cm]{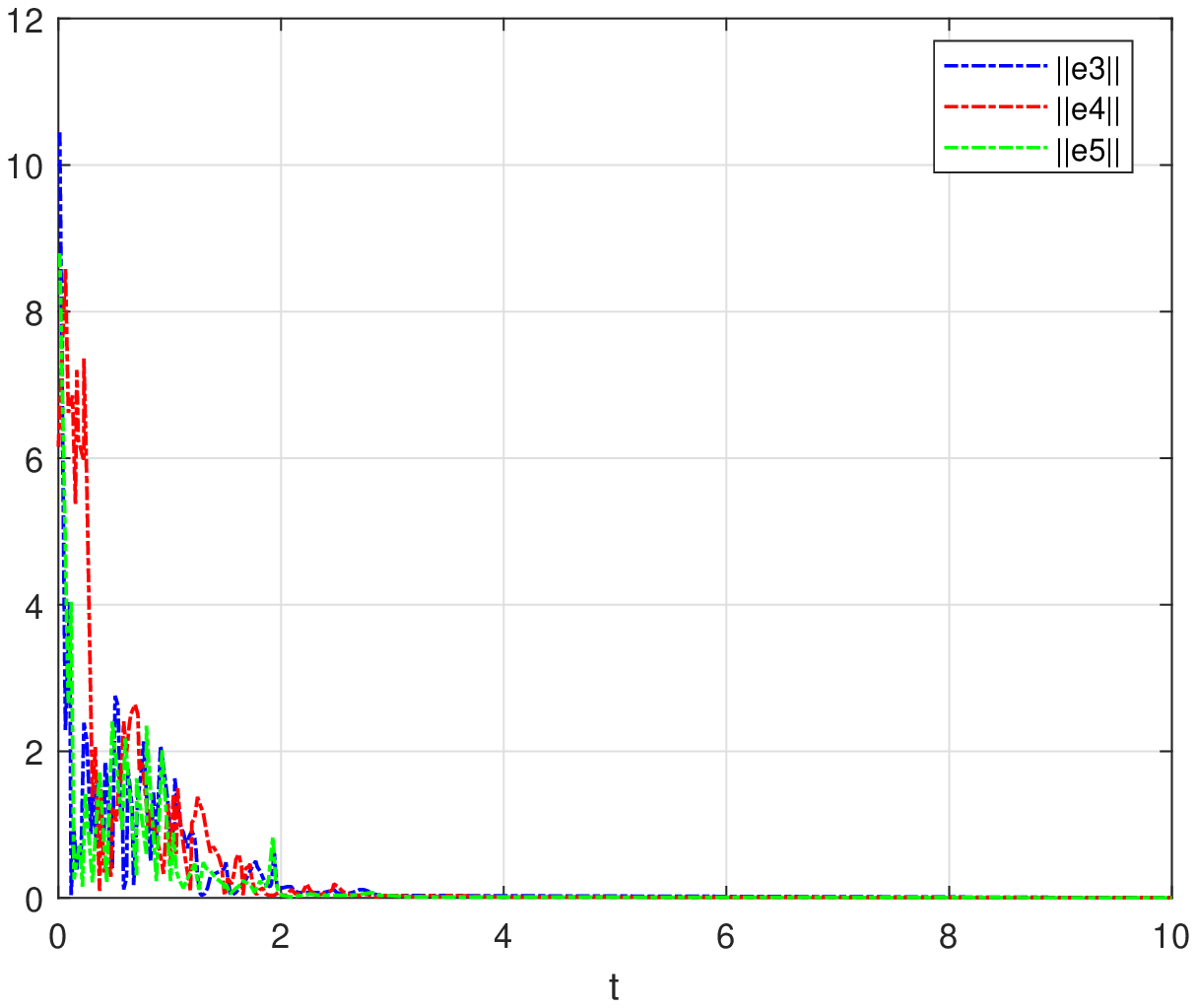}
\caption{The error dynamic trajectory when $i=3,4,5$ with the pinning impulsive controllers (6).}
\end{minipage}
\end{figure}

\begin{figure}[htbp]
\begin{minipage}[t]{0.45\linewidth}
\centering
\includegraphics[height=5cm,width=5cm]{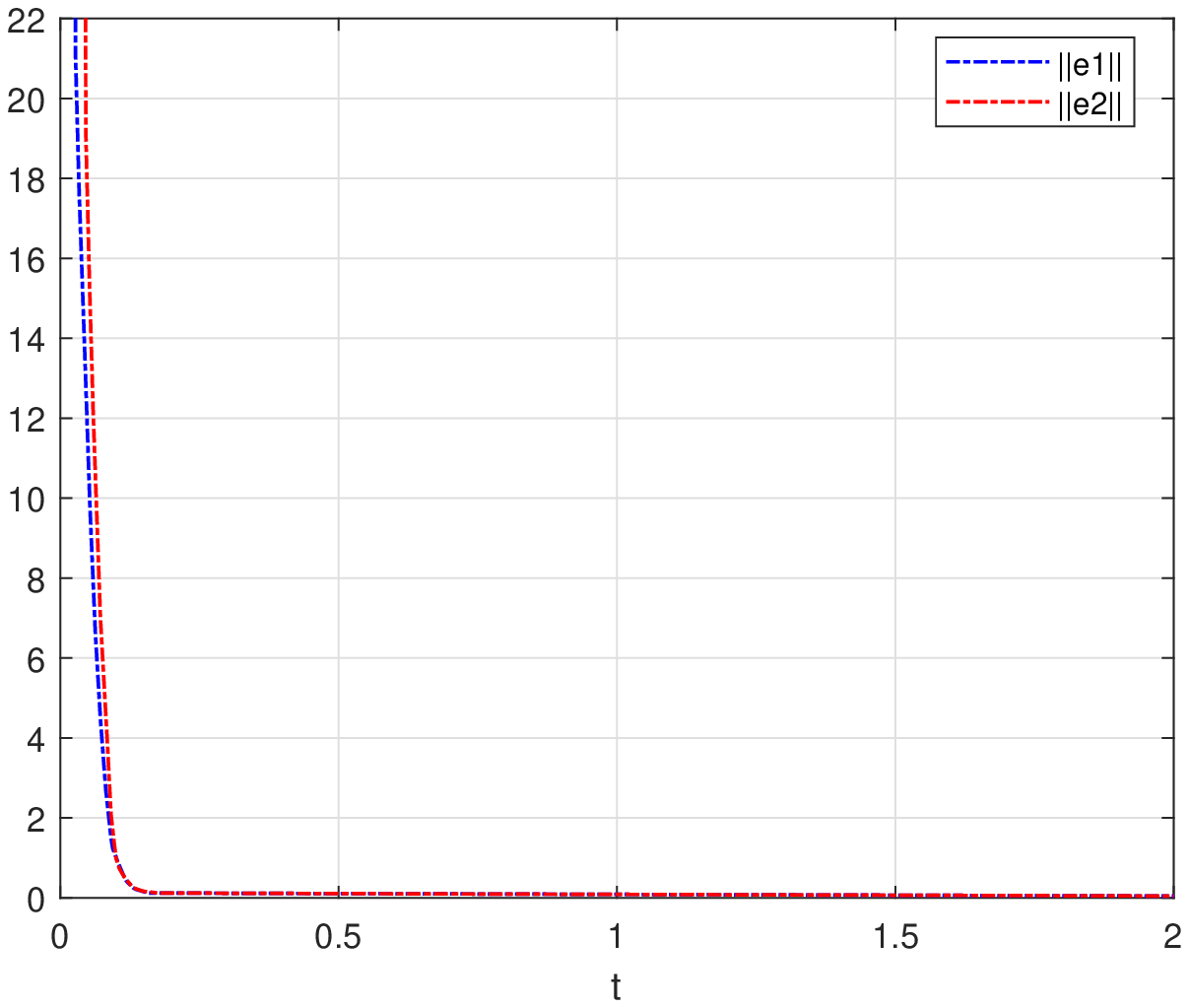}
\caption{The error dynamic trajectory $i=1,2$ with the hybrid controllers (23).}
\end{minipage}
\begin{minipage}[t]{0.45\linewidth}
\centering
\includegraphics[height=5cm,width=5cm]{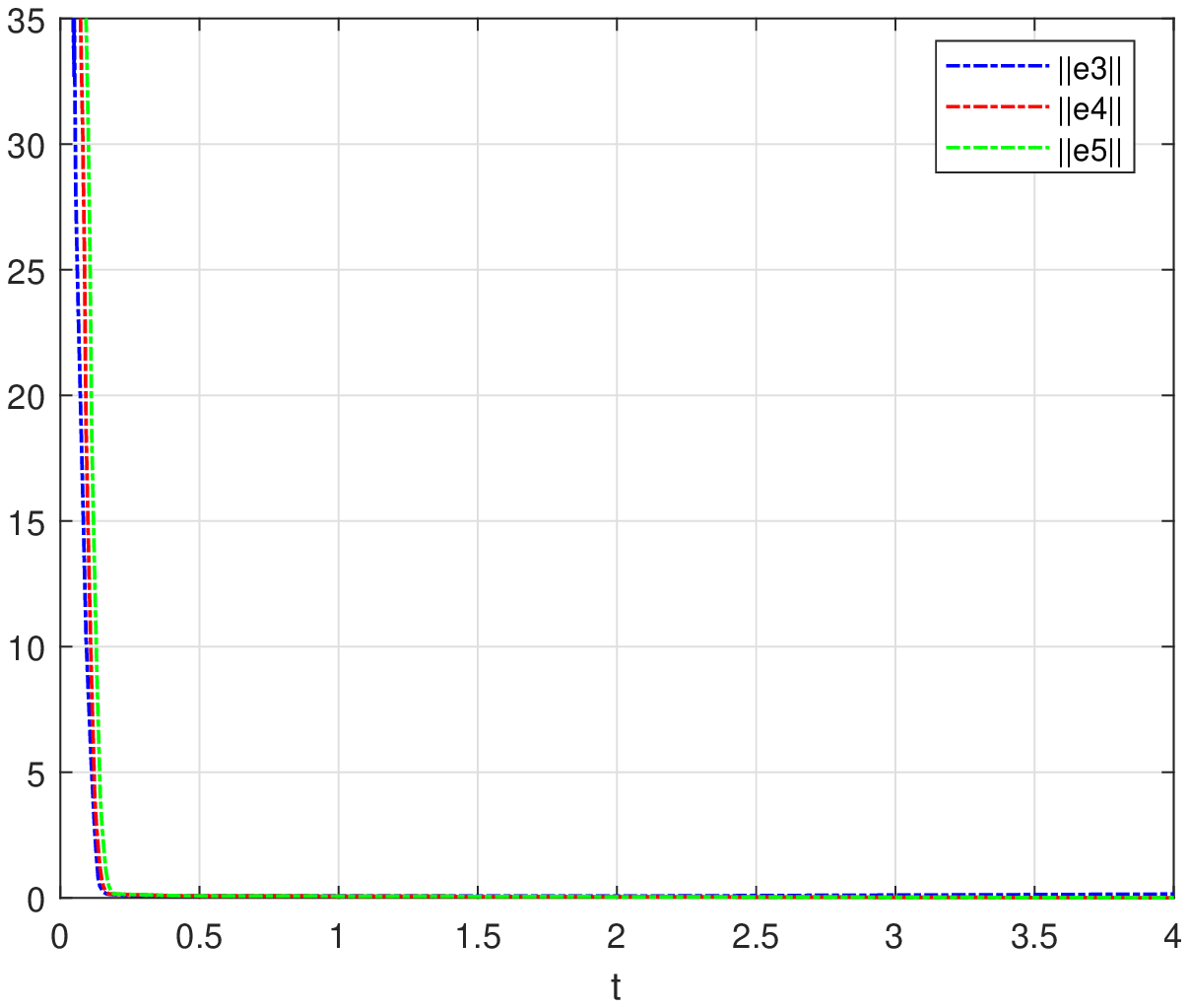}
\caption{The error dynamic trajectory when $i=3,4,5$ with the hybrid controllers (23).}
\end{minipage}
\end{figure}

\section{Conclusion}
\par In this article, the problem concerning asymptotic and finite-time cluster synchronization of coupled neural networks with time-varying delay via two different control schemes is mainly concerned. First, we select some typical nodes in the coupling neural networks, and apply the control of pinning impulse \cite{s33} which can achieve the goal of saving control costs and is helpful to obtain sufficient conditions for asymptotic cluster synchronization of the systems. Second, motivated by the goal of pinning control method which is operated by precedence controlling some key nodes, we introduce the hybrid controllers to deal with finite-time synchronization of the error system. Third, we show a numerical example by exhibiting three cases to show the correctness and effectiveness.
\par Using the properties of fractional-order calculus and the research technique of fractional-order energy function, the influence of two control schemes on the asymptotic and finite-time cluster synchronization of the fractional-order coupled neural networks with time-varying delays is studied. The result will emerge in the near future. What's more, complex-valued or quaternion-valued coupled neural networks with time-varying delays under the two control schemes in this paper have attracted the interest of many researchers. Hence, the future works also include asymptotic and finite-time cluster synchronization of complex-valued and quaternion-valued coupled neural networks with time-varying delays by pinning impulsive control and hybrid control.

\bibliographystyle{elsarticle-num}
\bibliography{text}

\begin{thebibliography}{10}

\bibitem{s1}
J.~Zha and C.~Li.
\newblock Synchronization of complex network based on the theory of
  gravitational field.
\newblock {\em Acta Physica Polonica B}, 50(1):87--114, 2019.

\bibitem{s2}
J.~Zhou, J.~Chen, J.~Lu, and J.~Lu.
\newblock On applicability of auxiliary system approach to detect generalized
  synchronization in complex network.
\newblock {\em IEEE Transactions on Automatic Control}, 62(7):3468--3473, 2017.

\bibitem{s3}
S.~Zhu, J.~Zhou, X.~Yu, and J.~Lu.
\newblock Bounded synchronization of heterogeneous complex dynamical networks:
  A unified approach.
\newblock {\em IEEE Transactions on Automatic Control}, 66(4):1756--1762, 2020.

\bibitem{s4}
Z.~Xu, P.~Shi, H.~Su, Z.~Wu, and T.~Huang.
\newblock Global h∞ pinning synchronization of complex networks with
  sampled-data communications.
\newblock {\em IEEE Transactions on Neural Networks and Learning Systems},
  29(5):1467--1476, 2018.

\bibitem{s5}
Y.~Wu and L.~Guo.
\newblock Enhancement of intercellular electrical synchronization by conductive
  materials in cardiac tissue engineering.
\newblock {\em IEEE Transactions on Biomedical Engineering}, 65(2):264--272,
  2018.

\bibitem{s6}
F.~Ren, F.~Cao, and J.~Cao.
\newblock Mittag-leffler stability and generalized mittag-leffler stability of
  fractional-order gene regulatory networks.
\newblock {\em Neurocomputing}, 160:185--190, 2015.

\bibitem{s7}
S.~Kumar, R.and~Sarkar, S.~Das, and J.~Cao.
\newblock Projective synchronization of delayed neural networks with mismatched
  parameters and impulsive effects.
\newblock {\em IEEE Transactions on Neural Networks and Learning Systems},
  31(4):1211--1221, 2020.

\bibitem{s8}
A.~Abdurahman and H.~Jiang.
\newblock Improved control schemes for projective synchronization of delayed
  neural networks with unmatched coefficients.
\newblock {\em International Journal of Pattern Recognition and Artificial
  Intelligence}, 34(3), 2020.

\bibitem{s9}
Adel O., X.~Wang, V.~Pham, and G.~Giuseppe.
\newblock Coexistence of identical synchronization, antiphase synchronization
  and inverse full state hybrid projective synchronization in different
  dimensional fractional-order chaotic systems.
\newblock {\em Advances in Difference Equations}, 2018(1):1--16, 2018.

\bibitem{s10}
Q.~Yang, H.~Wu, and J.~Cao.
\newblock Global cluster synchronization in finite time for complex dynamical
  networks with hybrid couplings via aperiodically intermittent control.
\newblock {\em Optimal Control Applications and Methods}, 41(4):1097--1117,
  2020.

\bibitem{s11}
L.~Liu, K.~Liu, H.~Xiang, and Q~Liu.
\newblock Cluster synchronization for directed complex dynamical networks via
  pinning control.
\newblock {\em Physica A: Statistical Mechanics and its Applications},
  545:123580, 2020.

\bibitem{s12}
A.~Hu, J.~Cao, M.~Hu, and L.~Guo.
\newblock Cluster synchronization of complex networks via event-triggered
  strategy under stochastic sampling.
\newblock {\em Complexity}, 434:99--110, 2015.

\bibitem{s13}
Y.~Wang, Z.~Ma, J.~Cao, and A.~Alsaedi.
\newblock Adaptive cluster synchronization in directed networks with
  nonidentical nonlinear dynamics.
\newblock {\em Complexity}, 21:380--387, 2016.

\bibitem{s14}
W.~Wu, W.~Zhou, and T.~Chen.
\newblock Cluster synchronization of linearly coupled complex networks under
  pinning control.
\newblock {\em IEEE Transactions on Circuits and Systems}, 56(4):829--839,
  2009.

\bibitem{s15}
Q.~Yang, H.~Wu, and J.~Cao.
\newblock Global cluster synchronization in finite time for complex dynamical
  networks with hybrid couplings via aperiodically intermittent control.
\newblock {\em Optimal Control Applications and Methods}, 41(4):1097--1117,
  2020.

\bibitem{s16}
J.~Cao, G.~Chen, and P.~Li.
\newblock Generalized analytical solutions and experimental confirmation of
  complete synchronization in a class of mutually-coupled simple nonlinear
  electronic circuits.
\newblock {\em Nonlinear Sciences}, 113:294--307, 2017.

\bibitem{s17}
A.~Koronovskii, O.~Moskalenko, V.~Ponomarenko, M.~Prokhorov, and A.~Hramov.
\newblock Binary generalized synchronization.
\newblock {\em Chaos Splitions and Fractals}, 83:133--139, 2016.

\bibitem{s18}
D.~Senthilkumar, M.~Lakshmanan, and J.~Kurths.
\newblock Phase synchronization in unidirectionally coupled ikeda time-delay
  systems.
\newblock {\em Nonlinear Sciences}, 164(11):35--44, 2008.

\bibitem{s19}
K.~Li, J.~Zhao, H.~Zhang, and X.~Li.
\newblock On successive lag synchronization of a dynamical network with delayed
  couplings.
\newblock {\em IEEE Transactions on Control of Network Systems}, page~1, 2021.

\bibitem{s20}
S.~Yang, C.~Hu, J.~Yu, and H.~Jiang.
\newblock Finite-time cluster synchronization in complex-variable networks with
  fractional-order and nonlinear coupling.
\newblock {\em Neural networks: the official journal of the International
  Neural Network Society}, 135:212--224, 2020.

\bibitem{s21}
S.~Liu, N.~Jiang, A.~Zhao, Y.~Zhang, and K.~Qiu.
\newblock Secure optical communication based on cluster chaos synchronization
  in semiconductor lasers network.
\newblock {\em IEEE Access}, 8:11872--11879, 2020.

\bibitem{s30}
L.~Pan, J.~Cao, U.~Al-Juboori, and M.~Abdel-Aty.
\newblock Cluster synchronization of stochastic neural networks with delay via
  pinning impulsive control.
\newblock {\em Neurocomputing}, 366(13):109--117, 2019.

\bibitem{s31}
P.~Liu, Z.~Zeng, and J.~Wang.
\newblock Asymptotic and finite-time cluster synchronization of coupled
  fractional-order neural networks with time delay.
\newblock {\em IEEE Transactions on Neural Networks and Learning Systems},
  PP(99), 2020.

\bibitem{s34}
J.~Mei, M.~Jiang, W.~Xu, and B.~Wang.
\newblock Finite-time synchronization control of complex dynamical networks
  with time delay.
\newblock {\em Communications in Nonlinear Science and Numerical Simulation},
  18(9):2462--2478, 2013.

\bibitem{s36}
P.~He, S.~Ma, and T.~Fan.
\newblock Finite-time mixed outer synchronization of complex networks with
  coupling time-varying delay.
\newblock {\em Chaos}, 22(4):175--423, 2012.

\bibitem{s35}
Z.~Hou, Jun P., Z.~Liu, J.~Zou, and J.~Luo.
\newblock Theory of functional differential equations.
\newblock {\em Journal of Changsha Railway University}, 2001.

\bibitem{s22}
X.~Jin, Z.~Wang, H.~Yang, Q.~Song, and M.~Xiao.
\newblock Synchronization of multiplex networks with stochastic perturbations
  via pinning adaptive control.
\newblock {\em Journal of the Franklin Institute}, 2021.

\bibitem{s23}
S.~Ding and Z.~Wang.
\newblock Synchronization of coupled neural networks via an event-dependent
  intermittent pinning control.
\newblock {\em IEEE Transactions on Systems, Man, and Cybernetics: Systems},
  pages 1--7, 2020.

\bibitem{s24}
W.~Li, J.~Zhou, J.~Li, T.~Xie, and J.~Lu.
\newblock Cluster synchronization of two-layer networks via aperiodically
  intermittent pinning control.
\newblock {\em IEEE Transactions on Circuits and Systems II: Express Briefs},
  page~1, 2020.

\bibitem{s25}
X.~Li, D.~Peng, and J.~Cao.
\newblock Lyapunov stability for impulsive systems via event-triggered
  impulsive control.
\newblock {\em IEEE Transactions on Automatic Control}, 65(11):4908--4913,
  2020.

\bibitem{s26}
X.~Li, X.~Yang, and J.~Cao.
\newblock Event-triggered impulsive control for nonlinear delay systems.
\newblock {\em Automatica}, 117, 2020.

\bibitem{s27}
W.~He, F.~Qian, and J.~Cao.
\newblock Pinning-controlled synchronization of delayed neural networks with
  distributed-delay coupling via impulsive control.
\newblock {\em Neural Networks}, 85:1--9, 2017.

\bibitem{s28}
X.~Xiong, X.~Yang, J.~Cao, and R.~Tang.
\newblock Finite-time control for a class of hybrid systems via quantized
  intermittent control.
\newblock {\em Sci. China Inf. Sci.}, 63(9):192201, 2020.

\bibitem{s29}
S.~Liu, Wu~H., and J.~Cao.
\newblock Fixed‐time synchronization for discontinuous delayed
  complex‐valued networks with semi‐markovian switching and hybrid
  couplings via adaptive control.
\newblock {\em International Journal of Adaptive Control and Signal
  Processing}, 34(10):1359--1382, 2020.

\bibitem{s37}
Z.~Wu, Q.~Ye, and D.~Liu.
\newblock Finite-time synchronization of dynamical networks coupled with
  complex-variable chaotic systems.
\newblock 2013.

\bibitem{s38}
H.~Wang, Z.~Han, Q.~Xie, and W.~Zhang.
\newblock Finite-time chaos control via nonsingular terminal sliding mode
  control.
\newblock {\em Communications in Nonlinear Science and Numerical Simulation},
  14(6):2728--2733, 2009.

\bibitem{s39}
L.~Liu and Q.~Liu.
\newblock Cluster synchronization in complex dynamical network of nonidentical
  nodes with delayed and non-delayed coupling via pinning control.
\newblock {\em Physica Scripta}, 94(4), 2019.

\bibitem{s40}
X.~Liu and T.~Chen.
\newblock Finite-time and fixed-time cluster synchronization with or without
  pinning control.
\newblock {\em IEEE Transactions on Cybernetics}, 48(99):240--252, 2018.

\bibitem{s32}
Q.~Wu, Z.~Jin, and X.~Lan.
\newblock Global exponential stability of impulsive differential equations with
  any time delays.
\newblock {\em Applied Mathematics Letters}, 23(2):143--147, 2010.

\bibitem{s33}
C.~Aouiti, M.~Bessifi, and X.~Li.
\newblock Finite-time and fixed-time synchronization of complex-valued
  recurrent neural networks with discontinuous activations and time-varying
  delays.
\newblock {\em Circuits Systems and Signal Processing}, pages 1--23, 2020.

\end{thebibliography}

\end{document}